# KERNEL ESTIMATION OF GREEK WEIGHTS BY PARAMETER RANDOMIZATION

By Romuald Elie, Jean-David Fermanian and Nizar Touzi

*CREST, BNP-PARIBAS, Ecole Polytechnique and Imperial College*

A Greek weight associated to a parameterized random variable $Z(\lambda)$ is a random variable $\pi$ such that $\nabla_\lambda E[\phi(Z(\lambda))] = E[\phi(Z(\lambda))\pi]$ for any function $\phi$. The importance of the set of Greek weights for the purpose of Monte Carlo simulations has been highlighted in the recent literature. Our main concern in this paper is to devise methods which produce the optimal weight, which is well known to be given by the score, in a general context where the density of $Z(\lambda)$ is not explicitly known. To do this, we randomize the parameter $\lambda$ by introducing an a priori distribution, and we use classical kernel estimation techniques in order to estimate the score function. By an integration by parts argument on the limit of this first kernel estimator, we define an alternative simpler kernel-based estimator which turns out to be closely related to the partial gradient of the kernel-based estimator of $\mathbb{E}[\phi(Z(\lambda))]$.

Similarly to the finite differences technique, and unlike the so-called Malliavin method, our estimators are biased, but their implementation does not require any advanced mathematical calculation. We provide an asymptotic analysis of the mean squared error of these estimators, as well as their asymptotic distributions. For a discontinuous payoff function, the kernel estimator outperforms the classical finite differences one in terms of the asymptotic rate of convergence. This result is confirmed by our numerical experiments.

**1. Introduction.** Let $\lambda$ be some given parameter in $\mathbb{R}^d$, and define the function
$$V^\phi(\lambda) := \mathbb{E}[\phi(Z(\lambda))],$$
where $Z(\cdot)$ is a parameterized random variable with values in $\mathbb{R}^n$ and $\phi: \mathbb{R}^n \to \mathbb{R}$ is a measurable function. In many applications we are interested in the nu-merical computation of the function $V^\phi(\lambda)$ for some parameter $\lambda^0$, together with the sensitivities of $V^\phi$ with respect to the parameter $\lambda$.

In particular, in the financial literature, $V^\phi$ represents the no-arbitrage price of a contingent claim, defined by the payoff $\phi(Z(\lambda))$, in the context of a complete market with prices measured in terms of the price of the nonrisky asset (so that the model is reduced to the zero-interest rate situation). The sensitivities of $V^\phi$ with respect to the parameter $\lambda$ are called *Greeks*, and are widely used by the practitioners in their *profit and loss* analysis. In the context of the Black–Scholes model, the derivative of the option price with respect to the current underlying asset price is the so-called *Delta*, and represents the number of shares of risky asset to be held at each time in order to realize a dynamic perfect hedge of the option. The *Gamma* is the second derivative of the option price, with respect to the underlying asset price. It is an indicator of the variation of the hedging portfolio. Another important *Greek* is the so-called *Vega* (although not a Greek letter!) which is the derivative of the option price with respect to the *volatility* coefficient (see, e.g., Hull [12] for more details).

We observe that the case of Bermudean options (i.e., American options with finite possible exercise times) is included in our framework by taking $\phi(Z(\lambda))$ as the value of the option at the first possible exercise time. The case of American options (with a continuous set of exercise times) can be covered by a limit argument, but requires a small time asymptotic analysis in order to control for the stability of the variance; see [8].

Given a numerical scheme for the computation of the function $V^\phi$, the first natural idea for the numerical computation of the *Greeks* is the finite differences approximation of the corresponding derivative. In addition to the generic standard error on the numerical computation of the expectation, this approximation leads to a biased estimator at a finite distance and appears to be inefficient for discontinuous payoff functions $\phi$. We refer to L'Ecuyer and Perron [7], Detemple, Garcia and Rindisbacher [5] or Milstein and Tretyakov [14] for a theoretical analysis of the rate of convergence of this estimator. Two direct methods for computing the *Greeks* have been presented by Broadie and Glasserman [3]: (i) the pathwise method, which consists in differentiating the random variable $\phi(Z(\lambda))$ inside the expectation operator, and (ii) the likelihood ratio method which reports the differentiation on the distribution of $Z(\lambda)$. The first method requires the computation of the gradient of the payoff function $\phi$, which is a serious limitation in practice as $\phi$ is typically highly complicated or even not differentiable; see also Giles and Glasserman [6] for further developments in this direction. As for the second method (ii), it was (apparently) restricted to the very special cases where the distribution of $Z(\lambda)$ is known explicitly. This difficulty was overcome by Fournié, Lasry, Lebuchoux, Lions and Touzi [9, 10] who exploited the Malliavin integration-by-parts formula to show that, for *smooth* random variables

$Z(\cdot)$,

$$\nabla_\lambda \mathbb{E}[\phi(Z(\lambda))] = \mathbb{E}[\phi(Z(\lambda))\pi], \tag{1.1}$$

where $\pi$, the so-called *Greek weight*, is a random variable independent of the pay-off function $\phi$. A quick overview of the notion of Greek weights is reported in Section 2. Further developments of the results of [9] were obtained by Gobet and Kohatsu-Higa [11]. The comparison of the above different methods is available in the survey paper of Kohatsu-Higa and Montero [13].

An important observation is that the set of Greek weights which satisfy (1.1) is a convex set of random variables. By an easy variance reduction argument, it is easily seen that the *score* $\pi^* := \nabla_\lambda \ln f(\lambda^0, Z(\lambda^0))$ minimizes $\mathbb{V}\mathrm{ar}[\phi(Z(\lambda))\pi]$, whenever the density $f(\lambda, z)$ of the random variable $Z(\lambda)$ exists and is sufficiently smooth. In general, the use of the Malliavin calculus does not lead to this optimal Greek weight, except in trivial cases where the density $f(\lambda, z)$ is explicitly known, which corresponds to the case covered by [3].

The main purpose of this paper is to focus on the use of the optimal Greek weight in order to estimate the corresponding Greek by the Monte Carlo method. To do this, our main idea is to randomize the parameter $\lambda$ and to rewrite $V^\phi$ as a regression function:

$$V^\phi(\lambda) := \mathbb{E}[\phi(Z(\Lambda))|\Lambda = \lambda],$$

where $Z(\Lambda)$ is a random variable with density $\varphi(\lambda, z) := \ell(\lambda^0 - \lambda)f(\lambda, z)$, and $\ell(\lambda^0 - \cdot)$ is some given randomizing distribution on the parameter $\lambda$ around $\lambda^0$. In other words, the random variable $Z(\Lambda)|\Lambda = \lambda$ has the same distribution as the random variable $Z(\lambda)$ defined by the density $f(\lambda, z)$. We next assume that our observations consist of a family $\{(\Lambda_i, Z_i), 1 \leq i \leq N\}$ of independent pairs $(\Lambda_i, Z_i)$ drawn in the density $\varphi$, and we define various kernel estimators of the Greek

$$\nabla_\lambda \mathbb{E}[\phi(Z(\lambda))]_{|\lambda=\lambda^0} = \mathbb{E}[\phi(Z(\lambda^0))s(\lambda^0, Z(\lambda^0))], \tag{1.2}$$

where $s(\lambda, z) := \nabla_\lambda \ln f(\lambda, z)$ is the score function. The first natural idea is to notice that

$$\mathbb{E}[\phi(Z(\lambda^0))s(\lambda^0, Z(\lambda^0))] = \mathbb{E}[\phi(Z(\Lambda))s(\Lambda, Z(\Lambda))|\Lambda = \lambda^0], \tag{1.3}$$

which is a usual regression function. Thus, a two-steps estimation method is proposed: we first perform a kernel-based estimator $\hat{s}$ of the score function, and then we define a kernel regression estimator of the Greek by substituting $\hat{s}$ to $s$. In the sequel the resulting estimator is referred to as the double kernel-based estimator and is denoted by $\tilde{\beta}$.

Our next kernel estimator of the Greek is based on a convenient integration-by-parts in (1.2). This leads to a much simpler estimator $\hat{\beta}$ which turns out to be closely related to the estimator $\check{\beta}$, obtained by direct differentiation of the classical kernel regression estimator of $V^\phi(\lambda) = \mathbb{E}[\phi(Z(\Lambda))|\Lambda = \lambda^0]$. These two estimators will be referred to as the single kernel-based estimators.

Let us observe that, unlike the so-called Malliavin Greek technique, our suggested estimators are biased but do not require any advanced mathematical calculation for their implementation. These two features are shared with the finite differences method. Also, intuitively, the randomization of the parameter $\lambda$ introduces an additional noise which may imply that our estimators are less accurate than their classical competitors. Our numerical results show indeed that the Malliavin Greek estimators are by far more accurate even in the case of an Asian option where the Greek weight is not optimal. Therefore, the main purpose of this paper is to provide a deep comparison between our estimators and the ones of finite differences.

Our three suggested estimators are defined precisely in Section 3, and their asymptotic properties are discussed in Section 4. We show that $\hat{\beta}$ and $\check{\beta}$ are asymptotically equivalent. The asymptotic properties of $\tilde{\beta}$ are derived under stronger conditions on the pay-off function $\phi$ and the kernel functions. The simultaneous choice of the bandwidth and the number of observations is also more restrictive in the latter case.

An important observation is that the two single kernel based estimators coincide if and only if the randomizing distribution $\ell$ is a truncated exponential distribution. In this case, by conveniently relating the support of the truncated exponential distribution to the kernel bandwidth, we observe that the rate of convergence is independent of the dimension of the parameter $\lambda$. We next solve the optimal choice of the randomizing distribution within this class by minimizing the corresponding mean square error.

Our asymptotic results imply the following main property of the single kernel based estimators: for a discontinuous payoff function $\phi$, the asymptotic rate of convergence of our estimator is better than the classical finite differences one, whenever the order of the kernel function is larger than some explicit threshold. In the case of a truncated exponential randomizing distribution, with support related to the kernel bandwidth, the single kernel based estimator has a better asymptotic rate of convergence whenever the order of the kernel function is larger than four.

Some numerical results are reported in Section 5. We estimate the delta of an European and an Asian digital call option. Our experiments show that the Malliavin-based estimators defined in [9] or [3] are the most efficient, as documented by the previous literature. As predicted by our theoretical asymptotic results, the single-kernel based estimator outperforms the finite differences one, but this is only observed for a large number of simulations. We believe that this does not restrict the interest in our new suggested method as this is just a matter of computer power, and the required number of simulations can be significantly reduced by using variance reduction techniques. For instance, the technique of antithetic variables applied to the randomizing density appears to be very efficient.

**2. The Greek weights set.** Throughout this paper we consider a complete probability space $(\Omega, \mathcal{F}, P)$. Let $Z(\lambda)$ be some random variable, valued in $\mathbb{R}^n$, depending on some finite-dimensional parameter $\lambda \in \mathbb{R}^d$, and set

$$V^\phi(\lambda) := E[\phi(Z(\lambda))] \qquad \text{for } \phi \in \mathbb{L}^\infty(\mathbb{R}^n, \mathbb{R}).$$

In order to simplify the presentation, we shall focus our attention on some fixed particular value $\lambda^0$ of $\lambda$, and we denote

$$Z^0 := Z(\lambda^0).$$

The chief goal of this paper is to devise efficient methods for the computation of the sensitivity parameter

$$\beta^0 := \nabla_\lambda V^\phi(\lambda^0),$$

for arbitrary functions $\phi$ chosen from a suitable large class. We assume that the distribution of $Z(\lambda)$ is absolutely continuous with respect to the Lebesgue measure, and we denote by $f(\lambda, z)$ the associated density, that is,

$$\mathbb{E}[\phi(Z(\lambda))] = \int \phi(z) f(\lambda, z)\, dz \qquad \text{for all } \phi \in \mathbb{L}^\infty(\mathbb{R}^n, \mathbb{R}).$$

Under mild smoothness assumptions on the density $f$, we directly compute that

$$\nabla_\lambda V^\phi(\lambda^0) := \frac{\partial V^\phi}{\partial \lambda}(\lambda^0) = E[\phi(Z^0) S^0], \qquad S^0 := s(\lambda^0, Z^0),$$

where the function $s$ is independent of $\phi$ and is explicitly given by

$$s(\lambda, z) := \nabla_\lambda \{\ln f(\lambda, z)\}.$$

This idea was introduced by Broadie and Glasserman [3] in the context of the Black–Scholes model where the density $f(\lambda, z)$ is explicitly known.

We shall always assume that

(2.1) $$\mathbb{E}|S^0|^2 < \infty.$$

Under this condition, the set

$$\mathcal{W} := \{\pi \in \mathbb{L}^2(\Omega, \mathbb{R}^d) : \nabla_\lambda V^\phi(\lambda^0) = E[\phi(Z^0)\pi] \text{ for all } \phi \in \mathbb{L}^\infty(\mathbb{R}^n, \mathbb{R})\}$$

is not empty. From the arbitrariness of $\phi \in \mathbb{L}^\infty(\mathbb{R}^n, \mathbb{R})$, it is immediately seen that

$$\mathcal{W} = \{\pi \in \mathbb{L}^2(\Omega, \mathbb{R}^d) : E[\pi | Z^0] = S^0\},$$

and therefore,

$$\begin{aligned}\mathbb{V}\mathrm{ar}[\phi(Z^0)\pi] &= E[\phi(Z^0)^2 E[\pi\pi'|Z^0]] - \nabla V^\phi(\lambda^0)\nabla V^\phi(\lambda^0)' \\ &\geq E[\phi(Z^0)^2 E[\pi|Z^0]E[\pi|Z^0]'] - \nabla V^\phi(\lambda^0)\nabla V^\phi(\lambda^0)' \\ &= E[\phi(Z^0)^2 S^0 S^{0'}] - \nabla V^\phi(\lambda^0)\nabla V^\phi(\lambda^0)' \\ &= \mathbb{V}\mathrm{ar}[\phi(Z^0)S^0].\end{aligned}$$

Hence,

$$S^0 \in \mathcal{W} \text{ is a minimizer of } \mathbb{V}\mathrm{ar}[\phi(Z^0)\pi], \qquad \pi \in \mathcal{W}.$$

Throughout this paper we call $S^0$ the optimal Greek weight. When the density function $f(\lambda, z)$ is not known, it was suggested in [9] to obtain (inefficient) Greek weights from the set $\mathcal{W}$ by exploiting the integration by-parts-formula from Malliavin calculus. Our main objective is to derive Monte Carlo estimators of the Greek value $\beta^0$, which asymptotically achieve the minimum variance, by using methods from nonparametric statistics to approximate the above optimal Greek weight $S^0$.

## 3. Kernel estimation and optimal Greek weight.

3.1. *Randomization of the parameter.* The main idea of this paper is to randomize the parameter $\lambda$ in order to estimate the Greek by the classical kernel estimation technique. This randomization can be exploited from two viewpoints. First, one can use it in order to estimate the optimal Greek weight, that is, the score function. An alternative viewpoint is to take advantage of the smoothness of the randomizing distribution in order to obtain an integration by parts formula similar to the Malliavin integration by parts technique. This technique is well known in the nonparametric statistics literature; see, for example, [1].

Let $\ell : \mathbb{R}^d \longrightarrow \mathbb{R}$ be some given probability density function, with support containing the origin in its interior, and set

$$\varphi(\lambda, z) := \ell(\lambda^0 - \lambda)f(\lambda, z) \qquad \text{for } \lambda \in \mathbb{R}^d \text{ and } z \in \mathbb{R}^n,$$

where $\lambda^0$ is the parameter of interest. We consider a sequence

(3.1) $(\Lambda_i, Z_i)_{1 \leq i \leq N}$ of $N$ independent r.v. with distribution $\varphi(\lambda, z)$,

so that, for any $i \leq N$, $\ell(\lambda^0 - \cdot)$ is the density of $\Lambda^i$ and $f(\Lambda^i, \cdot)$ is the conditional density of $Z^i$ given $\Lambda^i$.

REMARK 3.1. Notice that the simulation of $(\Lambda_i, Z_i)_{i \geq 1}$ can be performed easily even in cases where the density $\varphi$ can not be written explicitly. This applies typically to the case where $Z(\lambda) = X_T(\lambda)$, for some integer $T$, where $\{X_t(\lambda), t \in \mathbb{N}\}$ is a Markov chain with given transition density. Then, for a given value of $\lambda$, the simulation of $Z$ is easily feasible by usual methods. However, the marginal distribution of $Z(\lambda)$ is typically very complicated so that it is useless for the numerical computation of the score function $s(\lambda, z)$.

In this section we provide various estimation methods of $\beta$ based on nonparametric kernel methods. We then introduce the kernel function

$$K : \mathbb{R}^d \longrightarrow \mathbb{R} \quad \text{with} \int K = 1,$$

whose precise properties will be detailed at the beginning of Section 4.

3.2. *A first kernel estimator of the Greek.* The main idea is that the optimal weight $S^0$ requires a priori the knowledge of the probability density function $f(\lambda, z)$ and the associated *score* function $s(\lambda, z)$. Indeed, if these functions were explicitly known, then a natural nonparametric estimator of the Greek $\beta$ using the observations (3.1) is

$$(3.2) \qquad \bar{\beta}_N := \frac{1}{\ell(0) N h^d} \sum_{i=1}^{N} \phi(Z_i) s(\Lambda_i, Z_i) K\left(\frac{\lambda^0 - \Lambda_i}{h}\right).$$

Although $s$ is not explicitly known in our applications of interest, one could approximate it by means of an additional kernel estimator based on another kernel function $K$ defined on $\mathbb{R}^n$. We therefore introduce our first kernel-based estimator of $\beta$:

$$(3.3) \qquad \tilde{\beta}_N := \frac{1}{\ell(0) N h^d} \sum_{i=1}^{N} \phi(Z_i) \hat{s}_N^{-i}(\Lambda_i, Z_i) K\left(\frac{\lambda^0 - \Lambda_i}{h}\right),$$

where $s_N^{-i}$ is an approximation of $s$ given by

$$\begin{aligned}(3.4) \qquad \hat{s}_N^{-i}(\lambda, z) &:= \frac{\partial}{\partial \lambda} \ln \Bigg\{ \frac{1}{\ell(\lambda^0 - \lambda)(N-1) h^{d+n}} \\ &\qquad \times \sum_{j=1, j \neq i}^{N} K\left(\frac{\lambda - \Lambda_j}{h}\right) H\left(\frac{z - Z_j}{h}\right) \Bigg\} \\ &= \frac{\hat{\varphi}_\lambda^{-i}}{\hat{\varphi}^{-i}}(\lambda, z) + \frac{\nabla \ell(\lambda - \lambda^0)}{\ell(\lambda - \lambda^0)}\end{aligned}$$

and

$$(3.5) \qquad \hat{\varphi}^{-i}(\lambda, z) := \frac{h^{-d-n}}{N-1} \sum_{j=1, j \neq i}^{N} K\left(\frac{\lambda - \Lambda_j}{h}\right) H\left(\frac{z - Z_j}{h}\right),$$

$$
\hat{\varphi}_\lambda^{-i}(\lambda, z) := \nabla_\lambda \hat{\varphi}^{-i}(\lambda, z)
$$
(3.6)
$$
= \frac{h^{-d-n-1}}{N-1} \sum_{j=1, j \neq i}^{N} \nabla K\left(\frac{\lambda - \Lambda_j}{h}\right) H\left(\frac{z - Z_j}{h}\right).
$$

From a practical point of view, this estimator displays two drawbacks. First, its expression involves a product of two (possibly multidimensional) kernels $K$ and $H$. Thus, it suffers from the so-called "curse of dimensionality." Moreover, its calculation is time-consuming. In the subsequent subsections we introduce two alternative kernel estimators of $\beta$, which involve a single kernel function and a single summation.

From a theoretical point of view, we shall see that this estimator achieves the same rate of convergence as the two following ones but requires more stringent conditions, and involves heavy calculations.

3.3. *A simpler kernel estimator of the Greek.* We continue our discussion under the condition that

(3.7)  the kernel function $K$ has compact support.

We still consider the natural estimator given by (3.2). For fixed $h > 0$, it follows from the law of large numbers that

$$
\bar{\beta}_N \xrightarrow[N \to \infty]{} \mathbb{E}[\bar{\beta}_N]
$$
(3.8)
$$
= \frac{1}{\ell(0) h^d} \mathbb{E}\left[\phi(Z) s(\Lambda, Z) K\left(\frac{\lambda^0 - \Lambda}{h}\right)\right], \qquad P\text{-a.s.},
$$

where $(\Lambda, Z)$ is a random variable with distribution $\varphi(\lambda, z)$. Recalling the definition of $s$, and integrating by parts with respect to the variables $\lambda_1, \ldots, \lambda_d$, we see that, for $h > 0$ sufficiently small,

$$
\mathbb{E}[\bar{\beta}_N] = \frac{1}{\ell(0) h^d} \int \phi(z) K\left(\frac{\lambda^0 - \lambda}{h}\right) \ell(\lambda^0 - \lambda) \nabla_\lambda f(\lambda, z) \, d\lambda \, dz
$$
$$
= \frac{h^{-d-1}}{\ell(0)} \int \phi(z) \left(\nabla K\left(\frac{\lambda^0 - \lambda}{h}\right)\right.
$$
$$
\left. + h K\left(\frac{\lambda^0 - \lambda}{h}\right) \frac{\nabla \ell}{\ell}(\lambda^0 - \lambda)\right) \varphi(\lambda, z) \, d\lambda \, dz
$$
$$
= \frac{1}{\ell(0) h^{d+1}} \mathbb{E}\left[\phi(Z) \left(\nabla K\left(\frac{\lambda^0 - \Lambda}{h}\right)\right.\right.
$$
$$
\left.\left. + h K\left(\frac{\lambda^0 - \Lambda}{h}\right) \frac{\nabla \ell}{\ell}(\lambda^0 - \Lambda)\right)\right],
$$

where we used (3.7). This suggests the following simpler kernel estimator of the Greek $\beta$:

$$\hat{\beta}_N := \frac{1}{\ell(0)Nh^{d+1}} \sum_{i=1}^{N} \phi(Z_i)\Big(\nabla K\Big(\frac{\lambda^0 - \Lambda_i}{h}\Big) \quad\quad (3.9)$$
$$+ hK\Big(\frac{\lambda^0 - \Lambda_i}{h}\Big)\frac{\nabla \ell}{\ell}(\lambda^0 - \Lambda_i)\Big).$$

The asymptotic properties of $\hat{\beta}_N$ will be provided in Section 4.

3.4. *Differentiating the kernel estimator of the price.* We next start out from the natural kernel estimator of the price $V^\phi(\lambda)$:

$$\hat{V}_N^\phi(\lambda) := \frac{1}{Nh^d \ell(\lambda^0 - \lambda)} \sum_{i=1}^{N} \phi(Z_i) K\Big(\frac{\lambda - \Lambda_i}{h}\Big).$$

Differentiating $\hat{V}_N^\phi(\lambda)$ with respect to $\lambda$, we obtain our final kernel estimator of the Greek:

$$\check{\beta}_N := \frac{1}{\ell(0)Nh^{d+1}} \sum_{i=1}^{N} \phi(Z_i)\Big(\nabla K\Big(\frac{\lambda^0 - \Lambda_i}{h}\Big) \quad\quad (3.10)$$
$$+ hK\Big(\frac{\lambda^0 - \Lambda_i}{h}\Big)\frac{\nabla \ell}{\ell}(0)\Big).$$

Observe that our two estimators $\hat{\beta}_N$ and $\check{\beta}_N$ are closely related by

$$\check{\beta}_N = \hat{\beta}_N + \frac{1}{\ell(0)Nh^d} \sum_{i=1}^{N} \phi(Z_i) K\Big(\frac{\lambda^0 - \Lambda_i}{h}\Big)\Big(\frac{\nabla \ell}{\ell}(0) - \frac{\nabla \ell}{\ell}(\lambda^0 - \Lambda_i)\Big).$$

In particular,

$$\check{\beta}_N = \hat{\beta}_N \text{ whenever } \ell: l \mapsto e^{a_0 + a_1 \cdot l} \mathbf{1}_B(l) \quad\quad (3.11)$$

is a truncated exponential distribution,

for some parameters $a_0 \in \mathbb{R}$, $a_1 \in \mathbb{R}^d$ and some subset $B$ of $\mathbb{R}^d$ containing the origin in its interior.

The asymptotic properties of this third estimator will also be provided in Section 4.

**4. Asymptotic results.** We now compare the estimators defined in the previous section from the viewpoint of their asymptotic distributions. The heavy asymptotic analysis of the double kernel-based estimator is reported in [8], where it is shown to have the same asymptotic rate of convergence as the single kernel-based estimators, under more stringent conditions. Since the practical implementation of the double kernel-based estimator is, in addition, more time consuming, the discussion of this section will focus on the single kernel-based estimators.

We shall first show that the two single kernel-based estimators have equal asymptotic rates of convergence. Then, we derive the optimal choice of the number of simulations $N$ and the bandwidth $h$ of the kernel function $K$, by using the classical mean square error minimization criterion. We next specialize the discussion to the case of a truncated exponential randomizing distribution (3.11) with support defined by $B := [-\varepsilon, \varepsilon]^d$. In this setting, we observe that the rate of convergence of the kernel estimator is independent of the dimension of the parameter $\lambda$ for some convenient choice of $\varepsilon$ in terms of the bandwidth $h$. Finally, we discuss the optimal choice of the randomizing density $\ell$ within the class of truncated exponential distribution, and we provide a quasi-explicit characterization of the optimal truncated exponential randomizing distribution in the sense of the mean square error criterion.

Before stating our results, we recall that the order of the kernel function $K$ is defined as the smallest nonzero integer $p$ such that there exist some integers $(j_1, \ldots, j_p)$, $j_k \in \{1, \ldots, d\}$, such that

$$\int l_{\alpha_1} \ldots l_{\alpha_r} K(l) \, dl = 0 \quad \text{for } 0 < r < p, \alpha_k \in \{1, \ldots, d\}$$

and

$$\int l_{j_1} \ldots l_{j_p} K(l) \, dl \neq 0.$$

Typically, if $K$ is the product of $d$ even univariate kernels, then it is of order $p = 2$ (at least). The regularity hypothesis on the kernel function $K$ will be the following.

ASSUMPTION K. *The kernel function $K : \mathbb{R}^d \to \mathbb{R}$ is $C^1$, compactly supported, satisfies $\int K = 1$, and is of order $p \geq 2$.*

In the subsequent subsections we shall use the notation

$$(4.1) \quad \xi_K^p[\psi](\lambda, z) := \frac{(-1)^p}{p!} \sum_{j_1, \ldots, j_p = 1}^{d} \left( \int l_{j_1} \ldots l_{j_p} K(l) \, dl \right) \nabla^p_{\lambda_{j_1} \ldots \lambda_{j_p}} \psi(\lambda, z)$$

for every smooth function $\psi$ defined on $\mathbb{R}^d \times \mathbb{R}^n$. We shall also denote $A^\otimes := AA'$ for every matrix $A$.

4.1. *Asymptotic results for the single kernel-based estimators.* Our first result requires some regularity conditions on the density functions $f$ and $\ell$.

ASSUMPTION R1. For every $z$, the functions $f(\cdot, z)$ and $\ell$ are $p+1$ times differentiable, and for every integer $i \leq p$, the function $\lambda \longmapsto \nabla_\lambda^i \{\ell(\lambda^0 - \lambda)\nabla_\lambda f(\lambda, z)\}$ is continuous at $\lambda^0$ uniformly with respects to $z \in S$, for some subset $S$ s.t. $\mathrm{Supp}(\phi) \subset \mathrm{int}(S)$.

PROPOSITION 4.1. *Under Assumptions* K *and* R1, *as* $N \to \infty$ *and* $h \to 0$, *the bias and the variance of* $\hat{\beta}_N$ *satisfy*

$$\text{(4.2)} \qquad \mathbb{E}[\hat{\beta}_N] - \beta \sim C_1 h^p \quad and \quad \mathbb{V}\mathrm{ar}[\hat{\beta}_N] \sim \frac{\Sigma}{Nh^{d+2}},$$

*where*

$$\text{(4.3)} \qquad \begin{aligned} C_1 &:= \frac{1}{\ell(0)} \int \xi_K^p [\ell(\lambda^0 - \cdot) f_\lambda](\lambda^0, z) \phi(z)\, dz \quad and \\ \Sigma &:= \frac{\mathbb{E}[\phi^2(Z^0)]}{\ell(0)} \int \nabla K^\otimes. \end{aligned}$$

PROOF. By definition of $\hat{\beta}_N$, we have $\mathbb{E}[\hat{\beta}_N] = \mathbb{E}[\bar{\beta}_N]$. By (3.8), this provides

$$\psi(h) := \mathbb{E}[\hat{\beta}_N] = \frac{1}{\ell(0)h^d} \int \phi(z)\ell(\lambda^0 - \lambda)\nabla_\lambda f(\lambda, z) K\left(\frac{\lambda^0 - \lambda}{h}\right) d\lambda\, dz$$

$$= \frac{1}{\ell(0)} \int \phi(z)\ell(hl) f_\lambda(\lambda^0 - hl, z) K(l)\, dl\, dz.$$

Clearly, $\psi(0) = \int \phi(z) f_\lambda(\lambda^0, z)\, dz = \beta$. Moreover, since $K$ has compact support, it follows from Assumption R1 that the function $\psi$ is $p$ times differentiable at zero, with derivatives obtained by differentiating inside the integral sign, so that its $i$th iterated derivative denoted $\psi^{(i)}(0)$ are given by

$$\frac{(-1)^i}{\ell(0)} \sum_{j_1,\ldots,j_i=1}^d \left(\int l_{j_1}\ldots l_{j_i} K(l)\, dl\right)\left(\int \phi(z) [\nabla_{\lambda_{j_1},\ldots,\lambda_{j_i}}^i \{\ell(\lambda^0 - \cdot) f_\lambda\}](\lambda^0, z)\, dz\right)$$

for every $1 \leq i \leq p$. Since $p$ is the order of $K$, observe that $\psi^{(i)}(0) = 0$ for every $1 \leq i < p$, so that a Taylor expansion of $\psi$ provides the first part of the proposition.

As for the variance, we directly compute that

$$\mathbb{V}\mathrm{ar}[\hat{\beta}_N] = \frac{(v_1 - v_2^\otimes)}{Nh^{2d+2}\ell(0)^2},$$

where

$$v_1 := \mathbb{E}\bigg[\phi(Z)^2\bigg(\nabla K\bigg(\frac{\lambda^0 - \Lambda}{h}\bigg) + hK\bigg(\frac{\lambda^0 - \Lambda}{h}\bigg)\frac{\nabla \ell}{\ell}(\lambda^0 - \Lambda)\bigg)^{\otimes}\bigg],$$

$$v_2 := \mathbb{E}\bigg[\phi(Z)\bigg(\nabla K\bigg(\frac{\lambda^0 - \Lambda}{h}\bigg) + hK\bigg(\frac{\lambda^0 - \Lambda}{h}\bigg)\frac{\nabla \ell}{\ell}(\lambda^0 - \Lambda)\bigg)\bigg].$$

By a similar argument as in the first part of this proof, we compute that

$$v_1 = h^d \int \phi^2(z)\bigg(\nabla K(l) + hK(l)\frac{\nabla \ell}{\ell}(hl)\bigg)^{\otimes} \ell(hl) f(\lambda^0 - hl, z) \, dl \, dz$$

$$\sim h^d \ell(0)\bigg(\int \nabla K(l)^{\otimes} \, dl\bigg) \mathbb{E}[\phi^2(Z^0)].$$

The required result follows by observing that $v_2 = \mathrm{O}(h^{d+1})$. □

We are now ready for our first main result.

THEOREM 4.1. (i) *Let the conditions of Proposition* 4.1 *hold, and assume that*

(4.4) $$h \longrightarrow 0 \quad and \quad Nh^{d+2} \longrightarrow \infty \quad as \ N \to \infty.$$

*Then, with* $\Sigma$ *as in* (4.3), *we have* $\sqrt{Nh^{d+2}}(\hat{\beta}_N - \mathbb{E}[\hat{\beta}_N]) \longrightarrow \mathcal{N}(0, \Sigma)$ *in distribution.*

(ii) *In addition to the above conditions, assume that*

(4.5) $$Nh^{d+2+2p} \longrightarrow 0 \quad as \ N \to \infty.$$

*Then the bias vanishes and* $\sqrt{Nh^{d+2}}(\hat{\beta}_N - \beta) \longrightarrow \mathcal{N}(0, \Sigma)$ *in distribution.*

PROOF. We shall prove this result by verifying the Lyapounov conditions (see, e.g., Billingsley [2], page 44). Let $a$ be an arbitrary vector in $\mathbb{R}^d$, and define, for every $i = 1, \ldots, N$,

$$Y_i^N := \frac{1}{Nh^{d+1}\ell(0)}\phi(Z_i)\bigg(\nabla K\bigg(\frac{\lambda^0 - \Lambda_i}{h}\bigg) + hK\bigg(\frac{\lambda^0 - \Lambda_i}{h}\bigg)\frac{\nabla \ell}{\ell}(\lambda^0 - \Lambda_i)\bigg)$$

$$X_i^N := a'(Y_i^N - \mathbb{E}[Y_i^N]).$$

In view of Proposition 4.1, the only condition which remains to check in order to verify the Lyapounov conditions is the existence of $\delta > 2$ such that

(4.6) $$\sup_N \frac{1}{\sigma_N^\delta} \sum_{i=1}^N \mathbb{E}[|X_i^N|^\delta] < +\infty \quad \text{where } \sigma_N^2 := \mathbb{V}\mathrm{ar}\bigg[\sum_{i=1}^N X_i^N\bigg].$$

In order to prove (4.6), we start by observing from (4.2) that
$$\sigma_N^2 \sim \frac{1}{Nh^{d+2}\ell(0)} \mathbb{E}[\phi^2(Z^0)] \int |a'\nabla K(l)|^2 \, dl.$$

We next estimate by the Minkowski inequality and (4.2) that
$$\|X_i^N\|_\delta \leq \|a'Y_i^N\|_\delta + |a'\mathbb{E}[Y_i^N]|$$
$$= \|a'Y_i^N\|_\delta + \frac{1}{N}|a'\mathbb{E}[\hat{\beta}_N]|$$
$$\leq \frac{\sum_{i=1}^d \|\phi(Z)a_i(\nabla_i K((\lambda^0 - \Lambda)/h) + hK((\lambda^0 - \Lambda)/h)\nabla_i \ell/\ell(\lambda^0 - \Lambda))\|_\delta}{Nh^{d+1}\ell(0)}$$
$$+ \mathrm{O}\!\left(\frac{1}{N}\right)$$
$$\leq Const\left(\frac{h^{d/\delta}}{Nh^{d+1}} + \frac{1}{N}\right),$$

where the last estimate is obtained by a Taylor expansion with respect to the $h$ variable, in the neighborhood of the origin, following the method used in the proof of Proposition 4.1. Hence,
$$\frac{1}{\sigma_N^\delta}\sum_{i=1}^N \mathbb{E}[|X_i^N|^\delta] \leq Const\, N \frac{h^d}{(Nh^{d+1})^\delta}(Nh^{d+2})^{\delta/2} \leq \frac{Const}{(Nh^d)^{(\delta-2)/2}},$$

and condition (4.6) is satisfied when $Nh^d \to \infty$, as assumed in (4.4). Therefore, $\sum_{i=1}^N X_i^N$ is asymptotically Gaussian, with a variance matrix given by $\sigma_N^2$. By the arbitrariness of $a \in \mathbb{R}^d$, the required result follows from the Cramér–Wold device; see, for example, Theorem 25.5, page 405 in [4]. □

We next turn to the estimator $\check{\beta}$ which was defined as the gradient, with respect to $\lambda$, of the kernel based estimator $\hat{V}_N^\phi(\lambda)$ of the function $V_N^\phi(\lambda)$. The asymptotic properties of this estimator are obtained by following the techniques of the previous proofs and require the following regularity condition on the densities $f$ and $\ell$.

ASSUMPTION R2. For every $z$, the functions $f(\cdot, z)$ and $\ell$ are $p+1$ times differentiable, and for every integer $i \leq p$, the function $\lambda \longmapsto \nabla_\lambda^i \{\ell(\lambda^0 - \lambda)f(\lambda, z)\}$ is continuous at $\lambda^0$ uniformly with respect to $z \in S$, for some subset $S$ s.t. $\mathrm{Supp}(\phi) \subset \mathrm{int}(S)$.

PROPOSITION 4.2. *Under Assumptions* K *and* R2, *as* $N \to \infty$ *and* $h \to 0$, *the bias and the variance of* $\check{\beta}_N$ *satisfy*
$$\mathbb{E}[\check{\beta}_N] - \beta \sim C_2 h^p \quad \text{and} \quad \mathrm{Var}[\check{\beta}_N] \sim \frac{\Sigma}{Nh^{d+2}},$$

where $\Sigma$ is given by (4.3) and

$$C_2 := \frac{1}{\ell(0)} \int \left( \xi_K^p[\varphi_\lambda] + \frac{\nabla \ell}{\ell}(0) \xi_K^p[\varphi] \right)(\lambda^0, z)\phi(z)\, dz.$$

PROOF. The proof is essentially similar to that of Proposition 4.1. Recall that the estimators $\check{\beta}_N$ and $\hat{\beta}_N$ are related by

$$(4.7) \quad \check{\beta}_N = \hat{\beta}_N + \frac{1}{\ell(0)Nh^d} \sum_{i=1}^N \phi(Z_i) K\left(\frac{\lambda^0 - \Lambda_i}{h}\right) \left(\frac{\nabla \ell}{\ell}(0) - \frac{\nabla \ell}{\ell}(\lambda^0 - \Lambda_i)\right).$$

We start by analyzing the bias term. Recall from the proof of Proposition 4.1 that

$$\mathbb{E}[\hat{\beta}_N] = \frac{1}{\ell(0)} \int \phi(z) \ell(hl) f_\lambda(\lambda^0 - hl, z) K(l)\, dl\, dz.$$

We then deduce from (4.7) that

$$\mathbb{E}[\check{\beta}_N] = \frac{1}{\ell(0)} \int \phi(z) \left( \varphi_\lambda(\lambda^0 - hl, z) + \frac{\nabla \ell}{\ell}(0) \varphi(\lambda^0 - hl, z) \right) K(l)\, dl\, dz.$$

We now observe that Assumption R2 allows to derive an expansion in the $h$ variable of the above expression, near the origin, up to the order $p$. The coefficients of the expansion are obtained by simple differentiation inside the integral sign. Finally, since $p$ is the order of the kernel $K$, it is easily seen that the coefficients of $h^i$, $i < p$, in this expansion vanish, and the only nonzero coefficient is that of $h^p$.

The variance of $\check{\beta}_N$ is also treated by the same argument as in the proof of Proposition 4.1, and the dominant term in the expansion of the variance is easily seen to be the same as in that proof. $\square$

Proposition 4.2 says that $\hat{\beta}_N$ and $\check{\beta}_N$ have the same asymptotic variance, and the orders of their asymptotic biases are the same. Our next result states that these two estimators have exactly the same asymptotic distribution.

THEOREM 4.2. (i) *Let the conditions of Proposition* 4.2 *hold, and assume further that* (4.4) *holds. Then, with $\Sigma$ as in* (4.3), *we have* $\sqrt{Nh^{d+2}}(\check{\beta}_N - \mathbb{E}[\check{\beta}_N]) \longrightarrow \mathcal{N}(0, \Sigma)$ *in distribution.*

(ii) *Let* (4.5) *hold, in addition to the above conditions. Then the bias vanishes and*

$$\sqrt{Nh^{d+2}}(\check{\beta}_N - \beta) \underset{N \to \infty}{\longrightarrow} \mathcal{N}(0, \Sigma) \qquad \text{in distribution.}$$

PROOF. Define the sequence

$$Y_i^N := \frac{1}{Nh^{d+1}\ell(0)} \phi(Z_i) \left( \nabla K\left(\frac{\lambda^0 - \Lambda^i}{h}\right) + hK\left(\frac{\lambda^0 - \Lambda^i}{h}\right) \frac{\nabla \ell}{\ell}(0) \right),$$

and follow the lines of the proof of Theorem 4.1. $\square$

4.2. *Optimal choice of $N$ and $h$.* The two single kernel-based estimators $\hat{\beta}_N$ and $\check{\beta}_N$ have similar asymptotic properties. They both have a bias of order $h^p$, a variance of order $1/(Nh^{d+2})$ and a convergence in distribution at the rate $\sqrt{Nh^{d+2}}$. Therefore, the determination methods of the optimal $N$ and $h$ will be similar for both of them, and we only detail calculations for the estimator $\hat{\beta}_N$. Let the conditions of Proposition 4.1 hold. Then (4.2) holds, and we calculate an asymptotic equivalent for the mean square error between $\hat{\beta}_N$ and $\beta$:

$$\mathrm{MSE}(\hat{\beta}_N) := \mathbb{E}[|\hat{\beta}_N - \beta|^2] \sim \frac{\mathrm{Tr}(\Sigma)}{Nh^{d+2}} + h^{2p}|C_1|^2.$$

Minimizing the MSE in $h$, we get the asymptotically optimal bandwidth selector:

$$(4.8) \qquad \hat{h} = \left(\frac{(d+2)\mathrm{Tr}(\Sigma)}{2p|C_1|^2 N}\right)^{1/(d+2p+2)}.$$

Note that $\hat{h}$ is of order $N^{-1/(d+2p+2)}$, leading to an MSE of order $N^{-2p/(d+2p+2)}$. Similarly, the asymptotically optimal bandwidth selector for $\check{\beta}_N$ is

$$(4.9) \qquad \check{h} = \left(\frac{(d+2)\mathrm{Tr}(\Sigma)}{2p|C_2|^2 N}\right)^{1/(d+2p+2)}.$$

These results imply asymptotic theoretical choices for $h$ relative to $N$, but we may still encounter difficulties in the numerical calculation of $h$. Even if the optimal order of $h$ were known, we still need to evaluate the associated constant coefficients. From our empirical experiments, we observed that the accuracy of our estimators depends heavily on the choice of the bandwidth $h$, as usual in kernel estimation.

4.3. *The case of a uniform randomizing distribution.* We first study further the case where the randomizing density is uniform on the sphere of $\mathbb{R}^d$ centered at 0 with radius $\epsilon$:

$$\ell(l) \mapsto \frac{1}{(2\epsilon)^d}\mathbf{1}_{[-\epsilon,\epsilon]}(l).$$

Observe that this is a particular example from the truncated exponential distributions (3.11) for which the single kernel density estimators coincide:

$$\hat{\beta}_N = \check{\beta}_N = \frac{(2\epsilon)^d}{Nh^{d+1}}\sum_{i=1}^{N}\phi(Z_i)\nabla K\left(\frac{\lambda^0 - \Lambda_i}{h}\right).$$

Without loss of generality, we assume that the kernel $K$ has support on $[-1,1]^d$. We first rewrite Assumption R1 in the setting of this section.

ASSUMPTION R3. For every $z$, the function $f(\cdot, z)$ is $p+1$ times differentiable, and for every integer $i \leq p+1$, the function $\lambda \longmapsto \nabla^i_\lambda f(\lambda, z)$ is continuous at $\lambda^0$ uniformly with respect to $z \in S$, for some subset $S$ s.t. $\mathrm{Supp}(\phi) \subset \mathrm{int}(S)$.

PROPOSITION 4.3. *Let Assumptions K and R3 hold. Then, as $N \to \infty$, $h \to 0$ and $\epsilon \to 0$ with $\epsilon \geq h$, we have*

$$(4.10) \qquad \mathbb{E}[\hat{\beta}_N] - \beta \sim C_u h^p \quad and \quad \mathbb{V}\mathrm{ar}[\hat{\beta}_N] \sim N^{-1} h^{-d-2} \epsilon^d \Sigma_u,$$

*where*

$$(4.11) \qquad \begin{aligned} C_u &:= \int \xi_K^p[f_\lambda](\lambda^0, z) \phi(z) \, dz \quad and \\ \Sigma_u &:= 2^d \mathbb{E}[\phi^2(Z^0)] \int \nabla K^{\otimes}. \end{aligned}$$

PROOF. The proof is similar to the one of Proposition 4.1. Denoting by $\mathbf{1_d}$ the vector of $\mathbb{R}^d$ with unit component, we rewrite

$$\begin{aligned} \mathbb{E}[\hat{\beta}_N] &= \frac{1}{h^{d+1}} \int_{\mathbb{R}^n} \phi(z) \bigg( \int_{\lambda^0-\epsilon \mathbf{1_d}}^{\lambda^0+\epsilon \mathbf{1_d}} \nabla K\bigg(\frac{\lambda^0 - \lambda}{h}\bigg) f(\lambda, z) \, d\lambda \bigg) dz \\ &= \frac{1}{h} \int_{\mathbb{R}^n} \phi(z) \bigg( \int_{[-\epsilon/h, \epsilon/h]^d} \nabla K(u) f(\lambda^0 - uh, z) \, du \bigg) dz. \end{aligned}$$

Since $\epsilon \geq h$ and $K$ is supported on $[-1,1]^d$, we may replace in our last term the integration on $[-\frac{\epsilon}{h}, \frac{\epsilon}{h}]^d$ by an integration on $\mathbb{R}^d$, which is necessary to get the convergence of our estimator to $\beta^0$. Then, as in the proof of Proposition 4.1, an integration by parts followed by Taylor expansions gives us the expected equivalent of the bias. The same argument applies for the computation of the variance of $\hat{\beta}_N$. □

Sending $\epsilon$ to zero, we obtain the same asymptotic properties as in Proposition 4.1, as long as $\epsilon \geq h$. Therefore, the asymptotic optimal $\epsilon$ is simply the bandwidth $h$. The kernel-based estimator $\hat{\beta}_N^u$ associated with this optimal uniform density $\ell$ is then given by

$$(4.12) \qquad \hat{\beta}_N^u := \frac{2^d}{Nh} \sum_{i=1}^N \phi(Z_i) \nabla K\bigg(\frac{\lambda^0 - \Lambda_i}{h}\bigg),$$

and satisfies

$$(4.13) \qquad \mathbb{E}[\hat{\beta}_N^u] - \beta \sim C_u h^p \quad \text{and} \quad \mathbb{V}\mathrm{ar}[\hat{\beta}_N^u] \sim N^{-1} h^{-2} \Sigma_u,$$

with $C_u$ and $\Sigma_u$ defined in (4.11). Minimizing the corresponding mean square error, we obtain the optimal bandwidth

$$(4.14) \qquad h^u := \left(\frac{\operatorname{Tr}\Sigma_u}{p|C_u|^2 N}\right)^{1/(2p+2)}.$$

As in the study of the previous estimators, we also obtain a central limit theorem for the estimator $\hat{\beta}_N^u$.

THEOREM 4.3. (i) *Let the conditions of Proposition* 4.3 *hold in the particular case where $\epsilon = h$, and assume further that*

$$(4.15) \qquad h \longrightarrow 0 \quad and \quad Nh^2 \longrightarrow \infty \qquad as\ N \to \infty.$$

*Then, with $\Sigma_u$ as in* (4.11), *we have* $\sqrt{Nh^2}(\hat{\beta}_N^u - \mathbb{E}[\hat{\beta}_N^u]) \longrightarrow \mathcal{N}(0, \Sigma_u)$ *in distribution.*

(ii) *If, in addition, $Nh^{2p+2} \to 0$, then the bias vanishes and*

$$\sqrt{Nh^2}(\hat{\beta}_N^u - \beta) \longrightarrow \mathcal{N}(0, \Sigma_u) \qquad in\ distribution.$$

A remarkable feature of the above asymptotic result is that the rate of convergence is independent of the dimension $d$ of the parameter $\lambda^0$.

4.4. *The case of a truncated exponential randomizing distribution.* In this subsection we specialize the discussion to the one-dimensional case, and we consider a truncated exponential randomizing distribution:

$$\ell(l) := \theta \frac{e^{\theta l}}{e^{\theta \epsilon} - e^{-\theta \epsilon}} \mathbf{1}_{[-\epsilon,\epsilon]}(l),$$

with the parameter $\theta \in \mathbb{R}$, so that the two single kernel estimators associated to this density coincide:

$$\check{\beta}_N = \hat{\beta}_N = \frac{1}{\ell(0)Nh^{d+1}} \sum_{i=1}^{N} \phi(Z_i)\left(\nabla K\left(\frac{\lambda^0 - \Lambda_i}{h}\right) + \theta h K\left(\frac{\lambda^0 - \Lambda_i}{h}\right)\right).$$

Using the same line of arguments as in Proposition 4.3, we see that, under Assumptions K and R3, as $N \to \infty$, $h \to 0$ and $\epsilon \to 0$ with $\epsilon \geq h$, we have

$$(4.16) \qquad \mathbb{E}[\hat{\beta}_N] - \beta \sim C_e h^p \quad \text{and} \quad \mathbb{V}\mathrm{ar}[\hat{\beta}_N] \sim N^{-1}h^{-3}\epsilon \Sigma_e,$$

where $\Sigma_e := \Sigma_u$ defined in (4.11) and

$$C_e := \frac{(-1)^p}{p!}\left(\int u^p K(u)\,du\right)\sum_{k=1}^{p+1}\binom{p}{k-1}$$
$$(4.17)$$
$$\times \left(\int \nabla_\lambda^k f(\lambda^0, z)\phi(z)\,dz\right)(-\theta)^{p-k+1}.$$

Again, the asymptotic optimal $\epsilon$ is simply the bandwidth $h$ and the kernel-based estimator $\hat{\beta}_N^e$ associated with this optimal exponential density is given by

$$(4.18) \quad \hat{\beta}_N^e := \frac{e^{\theta h} - e^{-\theta h}}{\theta N h^2} \sum_{i=1}^N \phi(Z_i)\left(\nabla K\left(\frac{\lambda^0 - \Lambda_i}{h}\right) + \theta h K\left(\frac{\lambda^0 - \Lambda_i}{h}\right)\right).$$

The optimal bandwidth is obtained by minimizing the corresponding mean square error:

$$(4.19) \quad h^e := \left(\frac{\operatorname{Tr}\Sigma_e}{p|C_e|^2 N}\right)^{1/(2p+2)},$$

which leads to the following MSE:

$$(4.20) \quad \operatorname{MSE}(\hat{\beta}_N^e) = 2(p+1)p^{-p/(p+1)}[|C_e|^2(\operatorname{Tr}\Sigma_e)^p]^{1/(p+1)} N^{-p/(p+1)}.$$

As in Theorem 4.3, a central limit theorem for the estimator $\hat{\beta}_N^e$ can be derived.

REMARK 4.1. From the asymptotic viewpoint, the estimators based on the truncated exponential randomizing density differ by their bias, as the constants $C_e$ depends on $\theta$ while the variance $\Sigma_e = \Sigma_u$ is independent of $\theta$. The optimal truncated exponential randomizing density is then obtained by minimizing the squared bias, defined by the polynomial function $C_e^2$, with respect to $\theta$. In our numerical experiments of Section 5, this minimization is performed by classical Newton–Raphson iterations.

REMARK 4.2. Notice that, in both cases, the choice of the radius $\epsilon$ of $\ell$ depends on the kernel function $K$ only through its support. For instance, if $\operatorname{supp}(K) = [-M, M]^d$, then the optimal radius is $\epsilon = Mh$.

4.5. *Comparison with the finite differences estimators.* We first start by recalling the finite differences estimators. For ease of presentation, we let $d = 1$. The finite differences estimator of the parameter $\beta^0 := \nabla_\lambda \mathbb{E}[\phi(Z(\lambda^0))]$ is based on the finite differences approximation of the gradient

$$\nabla_\lambda \mathbb{E}[\phi(Z(\lambda^0))] \sim \frac{\mathbb{E}[\phi(Z(\lambda^0 + \alpha\varepsilon))] - \mathbb{E}[\phi(Z(\lambda^0 - (1-\alpha)\varepsilon))]}{\varepsilon},$$

where $\varepsilon > 0$ is a "small" parameter, and $\alpha \in [0, 1]$. The values $\alpha = 0$, 0.5 and 1 correspond respectively to the backward, centered and forward finite difference. The above finite difference approximation suggests the following finite differences estimator of $\beta$:

$$\hat{\beta}_N^{FD} = \frac{1}{N\varepsilon} \sum_{i=1}^N (\phi[Z^i(\lambda^0 + \alpha\varepsilon)] - \phi[Z^i(\lambda^0 - (1-\alpha)\varepsilon)]).$$

The asymptotic properties of these estimators were first studied by L'Ecuyer and Perron [7]. In the case where $\lambda \mapsto \phi[Z(\lambda)] \in C_b^3(\mathbb{R}^d)$, when $N \to \infty$ and $\varepsilon \to 0$ with $N^{1/4}\varepsilon \to 0$, they obtained a parametric rate of convergence:

$$\sqrt{N}(\hat{\beta}_N^{\text{FD}} - \beta) \xrightarrow[N\to\infty]{} \mathcal{N}(0, \Sigma_\alpha) \qquad \text{in distribution,} \qquad \text{for } \alpha = 0, \tfrac{1}{2} \text{ and } 1.$$

When the payoff function $\phi$ has a countable number of discontinuities, Detemple, Garcia and Rindisbacher [5] obtained the following central limit theorems:

For $\alpha = \tfrac{1}{2}$, when $N^{1/5}\varepsilon \to 0$,
$$N^{2/5}(\hat{\beta}_N^{\text{FD}} - \beta) \xrightarrow[N\to\infty]{} \mathcal{N}(0, \Sigma_\alpha) \qquad \text{in distribution.}$$

For $\alpha = 0, 1$, when $N^{1/3}\varepsilon \to 0$,
$$N^{1/3}(\hat{\beta}_N^{\text{FD}} - \beta) \xrightarrow[N\to\infty]{} \mathcal{N}(0, \Sigma_\alpha) \qquad \text{in distribution.}$$

In the general case $d \geq 1$, the finite differences estimators are defined componentwise, and therefore, the rate of convergence is not affected by the dimension $d$ of the parameter $\lambda^0$.

The main objective of this paragraph is to provide an asymptotic comparison of the single-kernel based estimator with the finite differences one. The key point of our single-kernel based estimators is that the differentiation with respect to the parameter $\lambda$ is reported on the density of $Z(\lambda)$ so that our asymptotic results do not involve the regularity of the pay-off function $\phi$. For any pay-off function $\phi$, and when $Nh^{d+2p+2} \longrightarrow 0$, we derived in Theorems 4.1 and 4.2 that

$$\sqrt{Nh^{d+2}}(\hat{\beta}_N - \beta) \xrightarrow[N\to\infty]{} \mathcal{N}(0, \Sigma) \qquad \text{in distribution,}$$

where $p$ is the order of the kernel function. Minimizing the corresponding MSE, we obtained in Section 4.2 an optimal $h$ of order $N^{-1/(d+2p+2)}$ which, of course, almost satisfies the condition required in the convergence in distribution. Therefore, taking a bandwidth $h$ of order $N^{-1/(d+2p+2)-2\delta/(d+2)}$ with $\delta > 0$ sufficiently small leads to a convergence in distribution at rate $N^r$ with $r := p/(d+2p+2) - \delta > 0$. Therefore, the single-kernel based estimators, with kernel of order $p > 2d + 4$ and $\delta$ sufficiently small, achieve a convergence rate of order $r > 2/5$. Hence, they outperform all the finite differences estimators in the case of discontinuous payoffs.

Notice that, by taking kernel functions of order $p$ sufficiently large, we can obtain a convergence rate in distribution as close as desired to the parametric rate $\sqrt{N}$.

REMARK 4.3. Consider the optimized kernel estimators $\hat{\beta}_n^u$ and $\hat{\beta}_n^e$, based on uniform or exponential density $\ell$ on the sphere with radius $h$, as derived in Section 4.3. Then, for $Nh^{2p+2} \to 0$, we obtain a rate of convergence of $\sqrt{Nh^2}$. Therefore, in order to outperform the finite differences estimators of a Greek associated to a discontinuous payoff function $\phi$, one just needs to use a kernel function of order $p > 4$.

**5. Numerical results.** In this section we present some numerical results obtained in the Black–Scholes model:

$$S_t^x := x \exp\left[\left(r - \frac{\sigma^2}{2}\right)t + \sigma W_t\right], \qquad t \geq 0, x > 0,$$

where $W$ is a standard Brownian motion on $(\Omega, \mathcal{F}, \mathbb{P})$ with values in $\mathbb{R}$, and $r \in \mathbb{R}$, $\sigma > 0$ are two given constants. We focus on the estimation of the so-called Delta:

$$\beta := \nabla_x \mathbb{E}[\phi(Z^x)],$$

where $Z^x = S_T^x$ for an European option and $Z^x = \int_0^T S_t^x \, dt$ for an Asian option. As in the previous sections, we denote by $f(x, \cdot)$ the density of $Z^x$.

We simulate independent observations $X_i$ distributed in the (optimal) exponential randomizing distribution $\ell$ on the sphere centered at $S_0 = x$ with radius $h$, as derived in Section 4.4. The single-kernel based estimator $\hat{\beta}_N^e$ is therefore given by (4.18).

5.1. *Computation of the optimal bandwidth.* As the "bumping" parameter $\epsilon$ for the finite differences estimator, the bandwidth in kernel estimation needs to be chosen carefully. The asymptotic results of Section 4 provide the expression of the asymptotic optimal bandwidth. For the truncated exponential randomizing distribution, we obtain

$$h^e = \left(\frac{\Sigma_e}{pC_e^2 N}\right)^{1/(2p+2)},$$

where $\Sigma_e = 2\mathbb{E}[\phi^2(Z^x)] \int (\nabla K)^2$ and

$$C_e := \frac{(-1)^p}{p!} \left(\int u^p K(u) \, du\right) \sum_{k=1}^{p+1} \binom{p}{k-1} \mathbb{E}\left[\phi(Z^x) \frac{\nabla_x^k f(x, Z^x)}{f(x, Z^x)}\right] (-\theta)^{p-k+1}.$$

Given a kernel function $K$, the coefficient $\Sigma_e$ can be estimated by a standard Monte Carlo procedure. We next focus on the estimation of the parameter

$$\mathcal{E}_k := \mathbb{E}\left[\phi(Z^x) \frac{\nabla_x^k f(x, Z^x)}{f(x, Z^x)}\right]$$

for a given $k \in \{1, \ldots, p+1\}$.

(i) Let $Z^x = S_T^x = xe^Y$, where $Y$ has a normal distribution with mean $m := (r - \frac{\sigma^2}{2})T$ and variance $\Sigma := \sigma^2 T$. Then it is easily checked that

$$\nabla_x^k f(x,z) = \left[\sum_{i=0}^k a_i^k d(x,z)^i\right] \frac{f(x,z)}{x^k},$$

where

(5.1) $$d(x,z) := \frac{\ln z - \ln x - m}{\Sigma},$$

and the coefficients $(a_i^j)_{(i,j)\in\{0,\ldots,k\}^2}$ are given by

(5.2) $$a_i^0 = \mathbf{1}_{\{i=0\}}, \qquad a_i^{j+1} = a_{i-1}^j - j a_i^j - \frac{i+1}{\Sigma} a_{i+1}^j,$$

with the convention $a_i^j = 0$ for $i < 0$ and $i > j$. Hence,

$$\mathcal{E}_k = \frac{1}{x^k} \mathbb{E}\left[\phi(Z^x)\left(\sum_{i=0}^k a_i^k d(x,Z^x)^i\right)\right],$$

and this parameter can be estimated by a straightforward Monte Carlo procedure.

(ii) In practice, the distribution function is unknown, and the calculation of the previous paragraph can not be used to estimate $\mathcal{E}_k$. We suggest to mimic the same principle as the usual Silverman's rule-of-thumb in kernel estimation (see Scott [15], e.g.): let $\hat{m}$ and $\hat{\Sigma}$ be two given estimates of the mean and variance $\ln(Z^x/x)$, respectively, and define $\hat{d}(x,z)$ and $(\hat{a}_i^j)_{(i,j)\in\{0,\ldots,k\}^2}$ by substituting $(\hat{m}, \hat{\Sigma})$ to $(m, \Sigma)$ in (5.1)–(5.2); then the coefficient $\mathcal{E}_k$ is approximated by

$$\hat{\mathcal{E}}_k = \frac{1}{x^p} \mathbb{E}\left[\phi(Z^x)\left(\sum_{i=0}^k \hat{a}_i^k \hat{d}(x,Z^x)^i\right)\right].$$

Once the coefficients $\mathcal{E}_k$ estimated for $1 \leq k \leq p+1$, the parameter $\theta$ is chosen through a numerical minimization; see Remark 4.1. In the particular case of an uniform randomizing distribution ($\theta = 0$), remark that only the estimation of $\mathcal{E}_{p+1}$ is necessary.

Therefore, the numerical procedure is divided in three steps: first, we estimate the terms detailed in the previous subsection $\Sigma_e$, $\mathcal{E}_k$, $\hat{m}$ and $\hat{\Sigma}_e$ through a Monte Carlo procedure with very few simulations. Then, we calibrate the parameter $\theta$ by minimization and we deduce the exponential optimal theoretical bandwidth. Finally, we estimate the delta of the option by means of a single-kernel based estimator with the estimated bandwidth.

REMARK 5.1. The numerical effort dedicated to the calculation of the optimal bandwidth parameter $h$ is also encountered in the classical finite differences method, as the optimal bumping parameter $\epsilon$ involves some a priori numerical simulations.

5.2. *Numerical comparison of the estimators.* We present here numerical results obtained for the estimation of the delta of an European and an Asian at-the-money digital calls, that is, with a payoff of the form $\phi(s) = \mathbf{1}_{s>K}$. Since this payoff function is discontinuous, the results of Section 4.5 show that the single-kernel based estimator achieves a better rate of convergence than the finite differences estimators, whenever the kernel has order $p > 4$. The main object of this section is to verify the empirical validity of these asymptotic results.

In order to compare their behavior, each estimator has been computed 200 times and their empirical distributions have been approximated by classical smoothing nonparametric estimation.

Our numerical experiments are performed with the following values of the parameters:

$$S_0 = 120, \quad r = 0, \quad \sigma = 0.2, \quad T = 1 \quad \text{and} \quad K = 120.$$

We use the following polynomial kernel functions of order 2, 4 and 6, respectively, with support on $[-1, 1]$:

$$K_2(u) = \tfrac{3}{4}(1 - u^2),$$
$$K_4(u) = \tfrac{15}{32}(1 - u^2)(3 - 7u^2),$$
$$K_6(u) = \tfrac{105}{256}(1 - u^2)(33u^4 - 30u^2 + 5).$$

From the viewpoint of computing time, kernel based or finite differences estimations with the same number of simulations are comparable. All the numerical tests have been realized in Visual C++ on a Pentium 4 xeon 3 GHz processor with 1 Gb of RAM.

*European digital call option.* In the context of the Black–Scholes model, it was observed by [9] that the optimal weight for European options can be obtained by means of the Malliavin integration by part formula, and coincides with the likelihood estimator introduced by [3]. Therefore, we are not hoping to compete with the Malliavin-based Monte Carlo estimator.

From our numerical experiments, we observed that the gain from using kernel estimators based on an exponential rather than a uniform randomizing distribution $\ell$ was very poor, especially when the order of the kernel function increases. From a numerical viewpoint, the gain obtained at most counter-balanced the numerical price of the minimization procedure. The examples presented here are therefore based on a uniform randomization distribution $\ell$.

The distributions of the different estimators based on $N = 10^6$ simulations are reported in Figure 1. The good performance of the Malliavin estimator is confirmed by our numerical experiments. However, we observe surprisingly that the three kernel based estimators are less accurate than the centered finite differences one, although their numerical computing times are comparable, of the order of 2 seconds. According to Section 4.5, the kernel of order 6 should perform better than the other ones, but this is not the case here. Actually, the terms $C_e$ and $\Sigma_e$ are such that the constant term of the mean square error increases very fast with the variability of K, which naturally increases with its order. For example, the MSE of the estimator based on the kernel of order 4 is ten times bigger than the one of the finite differences one, although they have the same rate of convergence. Furthermore, the optimal bandwidth $h$ increases with the order of the kernel, so that the asymptotic approximations become less accurate.

In order to further investigate this effect, we increase the number of simulations. Figure 2 shows the distribution of the finite differences estimator and the kernel based estimator of order 6 based on $N = 10^9$ simulations where each simulation takes approximately 30 minutes on our computer. In this case, we observe that the kernel based estimator of order 6 truly out-performs the finite differences one: its bias and its variance are two times smaller. This confirms the theoretical asymptotic results obtained in Section 4.5. We do not consider that the high number of simulations required is a serious restriction since it is just a matter of computer power or time given to the simulation. Furthermore, the good performance of the kernel based estimators of high order can be observed for a smaller number of simulations if we use in addition variance reduction technique. For example, by performing the antithetic variable technique with respect to the randomizing density $\ell$, we observe that the kernel based estimator of order 6 outperforms the finite differences estimator with $6 \times 10^7$ simulations, corresponding to a computer time of about 2 minutes.

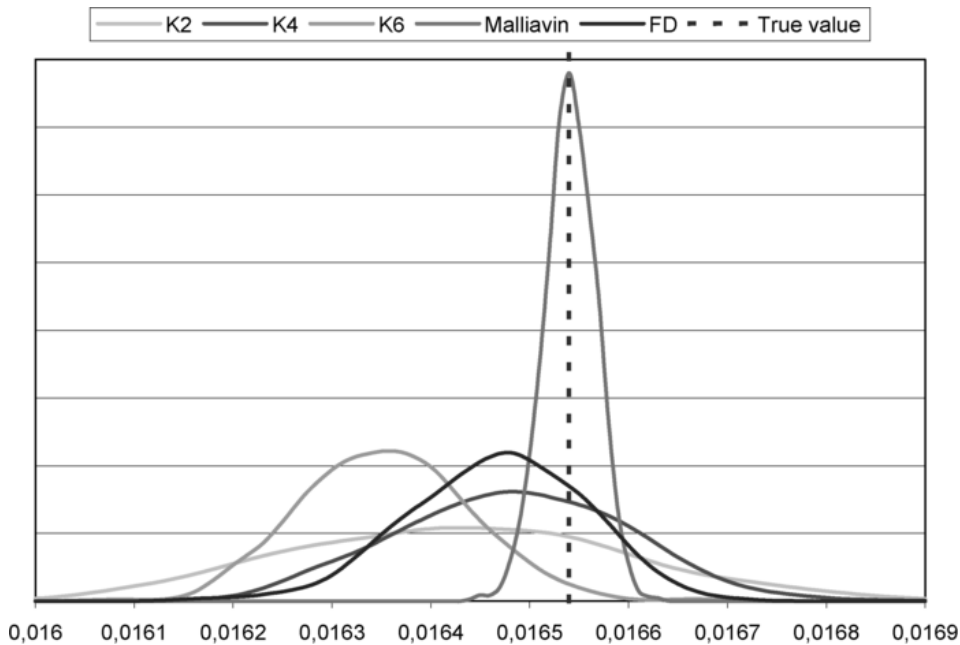

FIG. 1. *Delta of an European digital call, $N = 1$ million.*

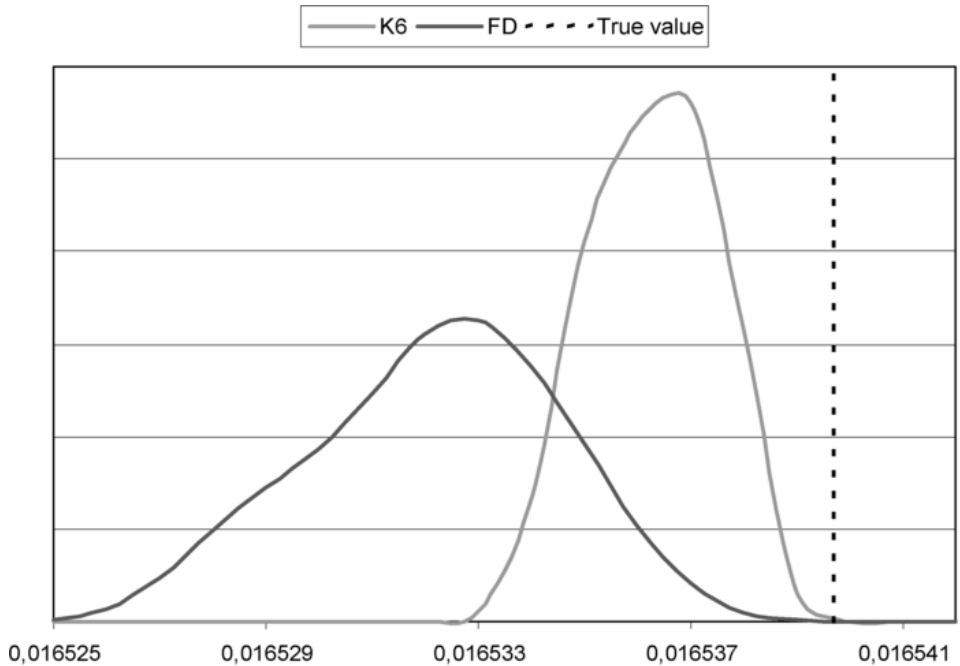

FIG. 2. *Delta of an European digital call, $N = 1$ billion.*

*Asian digital call option.* We next investigate the case of an Asian option, where the Malliavin integration by parts formula does not lead to the optimal weight; see [9]. The distribution of the different estimators based on $N = 10^6$ simulations are reported in Figure 3, where the "true value" of the Greek has been approximated by an unbiased Malliavin estimation with a very large number of simulations. Even if the Malliavin weight is not optimal, the Malliavin estimator still outperforms the other estimators. As for the European digital call, the finite differences estimator outperforms the kernel based estimators, but one simply requires more simulations in order to make the kernel estimator of order 6 more efficient than the finite differences one.

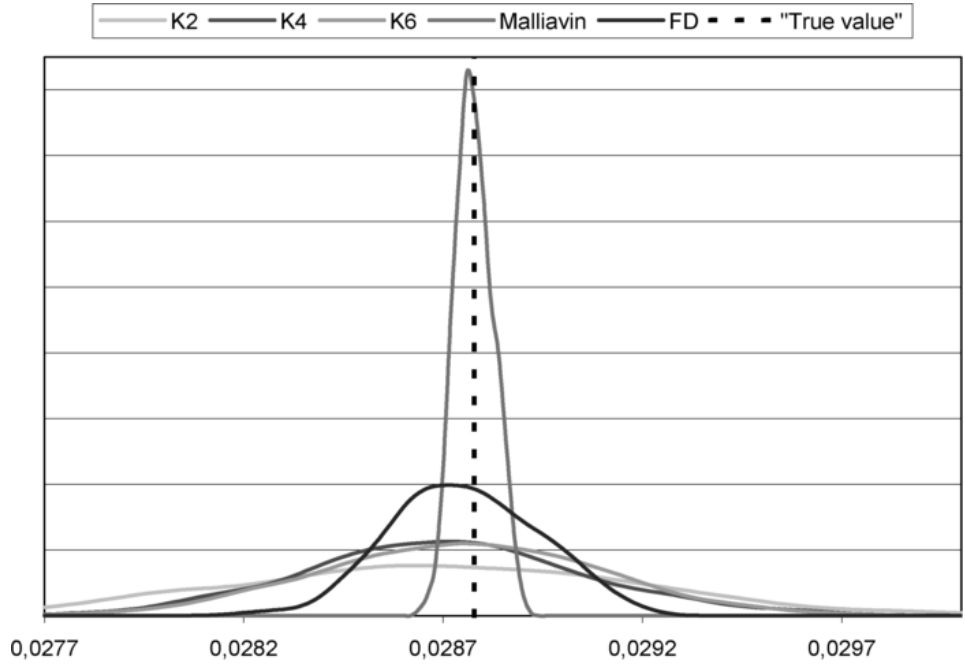

Fig. 3. *Delta of an Asian digital call, $N = 1$ million.*

*Conclusion* (*numerical results*). Other tests realized with different parameters, payoff functions or randomizing densities lead to rather similar results. Our kernel based estimator with order $p > 4$ of the delta of a digital option outperforms asymptotically the finite differences one, but one requires a large number of simulation to verify this fact empirically. Nevertheless, the high number of simulations required can be significantly reduced by means of variance reduction techniques. When the density of the underlying is unknown and the pay-off function is irregular, the Malliavin based estimator is still more efficient than the others. Nevertheless, in general, Malliavin weights are very difficult to derive analytically and this is precisely the advantage of the other estimators which are straightforward to implement.

## REFERENCES

[1] AIT-SAHALIA, Y. (1996). Nonparametric pricing of interest rate derivative securities. *Econometrica* **64** 527–560.
[2] BILLINGSLEY, P. (1968). *Convergence of Probability Measures.* Wiley, New York. MR0233396

[3] BROADIE, M. and GLASSERMAN, P. (1996). Estimating security prices using simulation. *Management Sci.* **42** 269–285.
[4] DAVIDSON, J. (1994). *Stochastic Limit Theory*. Oxford Univ. Press, New York. MR1430804
[5] DETEMPLE, J., GARCIA, R. and RINDISBACHER, M. (2005). Asymptotic properties of Monte Carlo estimators of derivatives. *Management Sci.* **51** 1657–1675.
[6] GILES, M. and GLASSERMAN, P. (2006). Smoking adjoints: Fast Monte Carlo Greeks. *Risk* 92–96.
[7] L'ECUYER, P. and PERRON, G. (1994). On the convergence rates of IPA and FDC derivative estimators. *Oper. Res.* **42** 643–656. MR1290513
[8] ELIE, R. (2006). *Double Kernel Estimation of Sensitivities*. To appear.
[9] FOURNIÉ, E., LASRY, J. M., LEBUCHOUX, J., LIONS, P. L. and TOUZI, N. (1999). Applications of Malliavin calculus to Monte Carlo methods in finance. *Finance and Stochastics* **3** 391–412. MR1842285
[10] FOURNIÉ, E., LASRY, J. M., LEBUCHOUX, J. and LIONS, P. L. (2000). Applications of Malliavin calculus to Monte Carlo methods in finance. II. *Finance and Stochastics* **5** 201–236. MR1841717
[11] GOBET, E. and KOHATSU-HIGA, A. (2003). Computation of Greeks for barrier and lookback options using Malliavin calculus. *Electron. Commun. Probab.* **8** 51–62. MR1987094
[12] HULL, J. (2002). *Options, Futures, and Other Derivatives.* Prentice Hall, Saddle River, New Jersey.
[13] KOHATSU-HIGA, A. and MONTERO, M. (2004). Malliavin calculus in finance. In *Handbook of Computational and Numerical Methods in Finance* (S. T. Rachev, ed.) 111–174. Birkhäuser, Boston. MR2083052
[14] MILSTEIN, G. N. and TRETYAKOV, M. V. (2005). Numerical analysis of Monte Carlo evaluation of Greeks by finite differences. *J. Comput. Finance* **8-3** 1–34.
[15] SCOTT, D. W. (1992). *Multivariate Density Estimation*. Wiley, New York. MR1191168

R. ELIE
CREST
10 RUE RUBENS
PARIS 75013
FRANCE
E-MAIL: elie@ensae.fr

J.-D. FERMANIAN
BNP-PARIBAS
10 HAREWOOD AVENUE
LONDON NW1 6AA
UNITED KINGDOM
E-MAIL: jean-david.fermanian@uk.bnpparibas.com

N. TOUZI
ECOLE POLYTECHNIQUE
UMR CNRS 7641
91128 PALAISEAU CEDEX
FRANCE
E-MAIL: touzi@ensae.fr